\input amstex
\input epsf
\documentstyle{amsppt}
\magnification=1200

\define\lra{\longrightarrow}
\redefine\qed{\hfill$\square$}
\define\diff{\text{Diff}}
\redefine\lim{\text{lim }}
\define\Th{\text{Th}}

\redefine\mod{\text{mod }}
\redefine\dot{\cdot}
\define\colim{\text{colim}\,}
\define\perb{\perp}
\define\nequiv{\not\equiv}
\redefine\coprod{\sqcup}
\define\bigcoprod{\bigsqcup}
\define\emb{\text{Emb}}
\define\den{{den}}

\define\FH{ H_{free}}
\define\C{\Bbb C}
\define\Z{\Bbb Z}
\define\F{\Bbb F}

\document

\heading
Divisibility of the stable Miller-Morita-Mumford classes  
\endheading

\vskip .2in
\author Soren Galatius, Ib Madsen, Ulrike Tillmann
\endauthor
\hfil {Soren Galatius, Ib Madsen, Ulrike Tillmann*} \hfil
\footnote[]{*The third author  was
supported by an Advanced Fellowship of the EPSRC.}
\hfil \break

\vskip .6in
\subheading {Abstract}
We 
determine the sublattice generated by the 
Miller-Morita-Mumford classes  $\kappa _i$ in the torsion free 
quotient of the integral cohomology ring of the stable mapping class group.
We further decide when the mod $p$ reductions $\kappa_i \in H^* (B\Gamma 
_\infty ; \Bbb F_p)$ vanish.

\vskip .2in
\subheading{Keywords}
Mapping class group, characteristic classes, surface bundles.

\vskip .2in
\subheading {MR classification} 57R20, 55P47.

\vskip .6in
\subheading{1. Introduction and results}

\vskip .1in
Let $\Gamma ^s_{g,b}$ denote the mapping class group of a surface of genus
$g$ with $b$ ordered boundary components and $s$  marked points.
We will supress $s$ or $b$ when their value is zero.
Gluing a   disk or a  torus with two boundary components to one of the boundary 
components induces homomorphisms
$$
\Gamma ^s_{g, b -1} \longleftarrow \Gamma ^s _{g,b} \lra \Gamma ^s_{g+1, b}.
\tag 1.1
$$
Recall that by Harer-Ivanov's stability theory  both homomorphisms induce 
a homology isomorphism in dimensions $*$ with $2* +1 < g$, cf. [H2], [I].
Let $\Gamma _{\infty} := \lim _{g\to \infty} \Gamma _{g,2}$ 
be the stable mapping class group.

\vskip .1in
Mumford in [Mu] introduced  certain tautological classes in the cohomology 
of moduli spaces of  Riemann surfaces. Miller [Mi]  and Morita [Mo] studied 
topological analogues: Let $e\in H^2( B \Gamma ^1_{g,b}; \Bbb Z)$ 
be the Euler class of the central extension
$$
\Bbb Z \lra \Gamma _{g,b+1} \lra \Gamma ^1_{g,b}
\tag 1.2 
$$
which is induced by gluing a disk with a marked point to one of the  boundary
components. 
Define 
$$
\kappa_i := \pi_ ! (e^{i+1}) \quad \in H^{2i} (B\Gamma _{g,b}; \Bbb Z)
$$
where $\pi_ !$ is the Umkehr (or integration along the fibre) map 
associated to the forgetful map $\Gamma _{g,b} ^1 \to \Gamma _{g,b}$.
These  correspond under the maps of (1.1) when $i>0$ and hence define classes
in $H^* (B\Gamma _\infty ; \Bbb Z)$. We will only be concerned with these 
stable classes in this paper.

\vskip .1in
By the proof of the Mumford conjecture [MW],
$$
H^*(B\Gamma _\infty; \Bbb Q) \simeq \Bbb Q [ \kappa_1, \kappa_2, \dots].
$$
In contrast, little is known about  $\kappa_i$ in integral cohomology though
it follows from [H1] that  $\kappa_1$ is
precisely divisible by $12$ (cf. [MT, p. 537]). 
We write 
$$
\FH ^* (B\Gamma_\infty) := H^* (B\Gamma _\infty; \Bbb Z)/ \text{Torsion}
$$
for the integral lattice in $H^*(B\Gamma _\infty; \Bbb Q)$.

\vskip .2in
\proclaim {Theorem 1.1}
Let $D_i$ be the maximal divisor  of $\kappa _i$ in $\FH ^* (B\Gamma _\infty)$.
Then for all $i\geq 1$
$$
D_{2i} = 2 \quad \quad \text{ and } \quad \quad D_{2i-1} = \den (\frac 
{B_i} {2i}). 
$$
\endproclaim
\vskip .1in
\noindent
Here $B_i$ denotes the $i$-th Bernoulli number and $\den$ is the function
that takes  a rational number
when expressed as a fraction in its lowest terms to its denominator.
It is well-known, cf.  [MSt; Appendix B], that
$\den (B_i)$ is the product of all primes $p$ such that $p-1$ divides
$2i$, and that a prime divides $\den (B_i/2i)$ if and only if it divides
$\den (B_i)$.
So in terms of
their $p$-adic valuation  the $D_i$ are determined by the formula
$$
\nu _p (D_i) = \left\{ \aligned
        1+\nu _p (i+1) \quad \quad &\text { if } i+1 \equiv 0 \quad \mod (p-1)\\
        0 \quad \quad &\text{ if } i +1  \nequiv 0 \quad \mod (p-1),
        \endaligned \right. 
\tag 1.3
$$
and $D_1 = 2^2\dot 3, $ $D_3= 2^3 \dot 3 \dot 5,$ $ D _5= 2^2\dot 3^2 \dot 7$
$  \dots$.

\vskip .2in

Our Theorem 1.1 is inspired by a conjecture of T. Akita [Ak] which we
also prove:

\vskip .2in
\proclaim{Theorem 1.2}
The element $\kappa_i$  in $H^{2i} (B\Gamma _\infty;
\Bbb F_p)$ 
vanishes if and only if
$i+1\equiv 0 \,  \text{mod } (p-1)$.
\endproclaim

\vskip .2in
\subheading{Remark 1.3}
The divisor $D_i$ of $\kappa_i$ in $\FH ^* (B\Gamma _\infty)$ is not 
necessarily equal  the  maximal divisor  $D_i^{\Bbb Z}$
 of $\kappa_i$ in integral cohomology
 $H^* (B\Gamma_\infty; \Bbb Z)$ but only provides an upper bound for it. 
However,   Theorem 1.2 gives
$$
p \text{ divides } D_i ^{\Bbb Z}\quad  \Longleftrightarrow \quad
p \text{ divides } D_i,
$$
which was strengthened by the first author in [G2] to
$$
p^2 \text{ divides } D_i ^{\Bbb Z} \quad  \Longleftrightarrow \quad
p^2 \text{ divides } D_i.
$$
It follows that for all even $i$ and for many odd $i$ 
($i= 1,5,9,13, \dots $), 
$D_i$ is indeed equal to $D_i^{\Bbb Z}$,  and one may expect that
$D_i = D_i ^{\Bbb Z}$  for
all $i \geq 1$.

\vskip .2in
\subheading{Remark 1.4} The integral lattice $H^*_{free}(B\Gamma_\infty)$
inherits a Hopf algebra structure.  The graded module of primitive
elements $P(H^*_{free}(B\Gamma_\infty))$ is a copy of $\Bbb  Z$ in each
even degree, and $\kappa_i$ is a primitive element of
$H^*_{free}(B\Gamma_\infty)$.  The structure of the Hopf algebra
$H^*_{free}(B\Gamma_\infty)$ is not completely understood at present,
but we have the following partial results.

\vskip .2in
\proclaim {Theorem 1.5} For
odd primes $p$ there is an isomorphism of Hopf algebras over the
$p$-local integers $\Bbb Z_{(p)}$
$$
  H^*_{free}(B\Gamma_\infty;\Bbb Z_{(p)}) \simeq H^*(BU;\Bbb Z_{(p)}).
$$
This fails for $p=2$.  In fact, the squaring map
$$
\xi: H^2_{free}(B\Gamma_\infty)\otimes \Bbb F_2 \to
      H^4_{free}(B\Gamma_\infty)\otimes \Bbb F_2
$$
is not injective, so the algebra $H^*_{free}(B\Gamma_\infty;\Bbb Z_{(2)})$
is not polynomial.
\endproclaim

\vskip .2in
In outline, the proofs of the  above theorems depend on previous 
results as follows.
For Theorem 1.2,
the proof of
the \lq\lq if" part in Section 3.3 
is a calculation of characteristic classes which  
relies on the fact that there is a map of infinite loop spaces
$\alpha : \Bbb Z \times B\Gamma ^+_\infty \to \Omega ^\infty 
\Bbb CP^\infty_{-1}$ 
(compare [T] and [MT] or Theorems 2.1 and 2.2 below).
The \lq \lq only if" part is implied by Theorem 1.1: if $p$ divides 
$\kappa_i$ then in particular it must divide its reduction to the
free part.

For Theorem 1.1,
we first establish a lower bound: $D_{2i} \geq 2$ by the \lq \lq if '' part of 
Theorem 1.2 and $D_{2i-1} \geq \den (B_i/2i)$
by a well-known relation between the $\kappa_i$
classes and the symplectic characteristic classes for surface bundles, 
(here stated as Theorem 4.2). 
The main theorem of [MT] provides an upper bound which is tight for $i$ even 
and precisely twice the lower bound for $i$ odd.
To eliminate the indeterminacy
of the factor 2, 
the main theorem of [MW] (Theorem 2.4 below), as well as calculations from
[G1] (Proposition 4.4)) and a stronger version of the main result of [MT]
(given in Theorem 2.2 and proved in Section 5) are used.

Theorem 1.5 is proved in Section 5.3.

\vskip .1in
Given the interest in the mapping class groups also
outside the topology community we have strived to make this
paper as self contained as possible. In particular we have spelled out some of
the more obscure parts of [MT].

\vskip .4in
\subheading{ 2. Spectrum cohomology and earlier results}

\vskip .2in
\subheading {2.1. Spectra and spectrum cohomology}
Let $E = \{ E_n, \epsilon_n\}$ be a CW-spec\-trum
\footnote{ If one does not assume the spaces to be CW-complexes then one
should assume that $\epsilon _n$ is a closed cofibration.}
in the sense
of [A1]: $E_n$ is a sequence of pointed CW-complexes
and $\epsilon_n: S E_n \to E_{n+1}$ a (pointed) isomorphism onto
a subcomplex, where $S(-)$ denotes suspension. The associated infinite 
loop space is the direct limit
$$
\Omega ^\infty E = \colim \, \Omega ^n E_n
$$
of the $n$-th loop space of $E_n$; the limit is taken over
the adjoint maps $\epsilon' _n: E_n \to \Omega E_{n+1}$.

The $k$-th homotopy group of $E$ is defined to be the direct limit of
$\pi_{n+k} (E_n)$. It is equal to the $k$-th homotopy group  of
the space $\Omega^\infty E$. In particular, the group of components of
$\Omega ^\infty E$ is the direct limit of 
$\pi_n (E_n)$. For $\alpha \in \pi_0 (\Omega ^\infty E)$ we 
let $\Omega ^\infty_\alpha E$ be the component determined by 
$\alpha$. In particular we write $\Omega ^\infty _0E$ for the component
of the zero element.

The homology and cohomology groups  of $E$ are 
$$
H^k(E) =\lim_{n\to \infty}
 \tilde H^{k+n} (E_n), \, \, H_k (E) =\colim_{n\to \infty} \tilde H_{k+n}(E_n)
$$
where the limits are induced from the maps $\epsilon _n$
together with the suspension isomorphisms. 
\footnote{ The $k$-th
spectrum cohomology  of $E$ is normally defined
as the group of
homotopy classes of degree $k$
spectrum maps from $E$ to the Eilenberg-MacLane
spectrum $K(\Bbb Z)$. This coincides with the above formula for $H^k(E)$
whenever for all $k$ the inverse system
$\lim _{n\to \infty}
\tilde H^{k+n} (E_n)$ satisfies the  Mittag-L\"offler condition
(see [A1]). In particular this will be the case when the structure
maps $\epsilon : S E_n \to E_{n+1} $ are $c(n)$-connected for some function
$c(n)$ with $c(n) \to \infty$ as $n\to \infty$.}
In contrast to homotopy groups
the cohomology groups of a spectrum are usually much simpler than the
cohomology groups of $\Omega ^\infty E$.

The evident evaluation map from $S^n\Omega ^n E_n$ to 
$E_n$ induces maps
$$
\sigma^*: H^* (E) \lra \tilde H^* (\Omega ^\infty_0 E), \, \, \, \sigma _*:
\tilde H_* (\Omega ^\infty _0E) \lra H_*(E).
\tag 2.1
$$
If we use field coefficients in the cohomology groups then
$H^* (\Omega ^\infty _0 E)$ is a connected Hopf algebra and 
the image of $\sigma ^*$ is contained in the graded vector space
$PH^* (\Omega ^\infty_0E)$ of primitive elements. We shall be particularly
concerned with the torsion free integral homology and cohomology groups
$$
\FH^*(\Omega ^\infty _0 E) = H^*(\Omega ^\infty _0E; \Bbb Z)/ \text{Torsion},
\quad 
H^{free}_*(\Omega ^\infty _0 E) = 
H_*(\Omega ^\infty _0E; \Bbb Z)/ \text{Torsion}. 
$$
They are lattices in $H^*(\Omega^\infty _0E; \Bbb Q)$ and $H_*(\Omega _0^\infty
E; \Bbb Q)$ and are dual Hopf algebras. Moreover, the image of 
$\sigma ^*$ is contained in the module of primitive elements
$$
\sigma ^*: \FH^*(E) \lra P(\FH^*(\Omega ^\infty_0E)),
$$
and dually $\sigma _*$ factors over the indecomposable  elements of
$H^{free}_* (\Omega ^\infty_0 E)$.

Given a pointed space $X$ we have the associated suspension
spectrum $S^\infty X$ whose $n$-th term is $S^nX$ with infinite loop space
$\Omega ^\infty S^\infty X$. 
There is an obvious inclusion
$i: X \to \Omega ^\infty S^\infty X$  inducing a splitting 
of $\sigma ^*$ (and $\sigma _*$):
$$
H^*(S^\infty X ) \overset {\sigma^*} \to \lra \tilde
 H^* (\Omega ^\infty_0 S^\infty
X) \overset {i^*} \to \lra \tilde H^* (X)
\tag 2.2
$$
is the suspension isomorphism.
$\Omega ^\infty S^\infty X$ is the free infinite loop space on $X$
and satisfies the universal property that any pointed
map from $X$ to some infinite     
loop space $Y$ can be extended in a unique way up to homotopy
to a map of infinite loop spaces from 
$\Omega ^\infty S^\infty X$ to $Y$. 

The spectra of most relevance to us are $\Bbb CP^\infty _{-1}$ and the 
suspension spectrum $S^\infty \Bbb CP^\infty _{+} $ of
$\Bbb CP^\infty \coprod \{ +\}$.
We recall the definition of the former. There are two 
complex vector bundles over the complex
projective $n$-space $\Bbb CP^n$, namely the tautological line bundle
$L_n$ and its $n$-dimensional complement $L_n^\perb$ in 
$\Bbb CP^n \times \Bbb C^{n+1}$. Its Thom space (or  one point compactification)
is denoted by $\Th (L_n^\perb)$. Since the restriction of $L_n^\perb$
to $\Bbb CP^{n-1} \subset \Bbb CP^n$ is equal to $L_{n-1}^\perb \oplus
\Bbb C$ where $\Bbb C$ denotes the trivial line bundle over
$\Bbb CP^{n-1}$ we get a map 
$$
\epsilon: S^2 \Th (L_{n-1} ^\perb) \lra \Th (L_n^\perb).
$$
The spectrum $\Bbb CP^\infty _{-1} $ has
$$
(\Bbb CP^\infty_{-1} )_{2n} = \Th (L_{n-1}^\perb), \, \, \, \quad 
(\Bbb CP ^\infty _{-1}) _{2n+1} = S \Th (L_{n-1} ^\perb)
$$
and the structure map $\epsilon _{2n+1}$ is given by the above $\epsilon$.
The associated infinite loop space is
$$
\Omega ^\infty\Bbb CP^\infty _{-1}
= \colim _{n \to \infty} \, \,  \Omega ^{2n}\Th (L_{n-1}^\perb).
$$
The inclusion of $L_{n-1} ^\perb $ into $L_{n-1} ^\perb \oplus L_{n-1} 
= \Bbb CP^{n-1} \times \Bbb C^n$
via the zero section of $L_{n-1}$
induces a map from $\Th (L_{n-1} ^\perb)$ into $S^{2n} (\Bbb CP^{n-1}_+)$
and hence a map
$$
\omega: \Omega ^\infty \Bbb CP ^\infty _{-1} \lra \Omega ^\infty S^\infty
(\Bbb CP^\infty _+).
$$
This map fits into a fibration sequence
$$
\Omega ^\infty \Bbb CP^\infty _{-1} \overset \omega \to \lra
\Omega ^\infty S^\infty (\Bbb CP^\infty _+) \overset \partial \to \lra
\Omega ^\infty S^{\infty -1}
\tag 2.3
$$
where the right-hand term is the direct limit of $\Omega ^nS^{n-1}$ [R].
Indeed the inclusion of a fibre $\Bbb C^n \to L_n ^\perb$
induces a map $S^{2n} \to \Th (L_n^\perb)$ and gives rise to a cofibre 
sequence of spectra $S^\infty (S^{-2}) \to \Bbb C P ^\infty _{-1} 
\to S^\infty \Bbb CP^\infty _+ \to S^\infty (S^{-1})$.
(2.3) is the associated fibration sequence of infinite loop spaces.

The component groups of (2.3) are
$$
0 \lra \Bbb Z \overset {\pi_0(\omega)} \to \lra \Bbb Z
\overset {\pi _0 (\partial)} \to \lra \Bbb Z/2 \lra 0
$$
so $\pi_0 (\omega) $ is multiplication by $\pm 2$, depending on
the choice of generators.  There is a canonical splitting of infinite 
loop spaces
$$
\Omega ^\infty S^\infty (\Bbb CP^\infty _+) \simeq \Omega ^\infty S^\infty
(\Bbb CP^\infty) \times \Omega ^\infty S^\infty.
\tag 2.4
$$
We fix the generator of $\pi_0 \Omega ^\infty S^\infty (\Bbb C P^\infty _+)$
to be the element that  maps  to +1 under the isomorphisms
$$
\pi_0( \Omega ^\infty S^\infty (\Bbb CP^\infty _+)) \overset {\pi_0 (c)}
\to \lra \pi_0 (\Omega ^\infty S^\infty) \overset {\text{degree}} \to \lra
\Bbb Z,
$$
where $c$ collapses $\Bbb CP^\infty$ to the non-base point of $S^0$.
We fix  the generator of $\pi_0 (\Omega ^\infty \Bbb C P^\infty _{-1})$ so that 
$\pi_0 (\omega) $ is multiplication by $-2$.

\vskip .3in
\subheading{ 2.2. Review of  results used}
Our divisibility result of Theorem 1.1 is based upon
the following three theorems.

\vskip .2in
\proclaim{Theorem 2.1. [T]} The spaces $\Bbb Z \times B\Gamma^+_\infty$ and
$B\Gamma ^+_\infty (= \{0\} \times B\Gamma ^+_\infty)$ are infinite loop spaces.
\endproclaim

Here the superscript (+) denotes Quillen's plus construction, cf. [B].
The product structure can be described as follows:
We may view $\Gamma _{g,2}$ as the mapping class group of
surfaces with one incoming and one outgoing boundary component. Gluing the
incoming boundary component of one surface to the
outgoing component of the other
defines a map
$$
\Gamma _{g,2} \times \Gamma _{h,2} \lra \Gamma _{g+h,2}
$$
and a corresponding map of classifying spaces that makes the disjoint union
$\bigcoprod B\Gamma _{g,2} $ over all $g\geq 0$ into a topological monoid.
Consider the map
$$
\bigcoprod _{g\geq 0} \, B\Gamma _{g,2} \lra \Bbb Z\times B\Gamma ^+_\infty
$$
that sends $B\Gamma _{g,2}$ into the component $\{g\} \times B\Gamma ^+_\infty$
by the stabilization map (1.1) followed by the map into the plus construction.
The infinite loop space structure on $\Bbb Z\times B\Gamma ^+_\infty$ is
compatible with the monoidal structure on $\bigcoprod B\Gamma _{g,2}$, and  the
induced map
$$
\Omega B (\, \bigcoprod _{g\geq 0}
B\Gamma _{g,2} ) \lra \Bbb Z\times B\Gamma ^+_\infty
$$
is a homotopy equivalence of loop spaces.
We refer to [T] for details.

\vskip .1in
To state the next result,
for each prime  $p$ we pick a positive integer $k= k(p)$
so that $-k$ reduces to a generator of the units $(\Bbb Z /p^2) ^\times$ when
$p$ is odd. We pick $k=5$ when $p=2$. Write $\psi ^{-k}$ for the self map 
of $\Bbb CP^\infty$ that multiplies by $-k$ on the second cohomology
group.  Composing with the inclusion into $\Omega ^\infty S^\infty
(\Bbb CP^\infty)$ and using the loop sum we have  a map
$$
1+ k\psi ^{-k} : \Bbb CP^\infty \lra \Omega ^\infty S^\infty (\Bbb CP^\infty),
$$
and, using the universal  property of free infinite loop spaces,
 a unique  extension to a self map
of $\Omega ^\infty S^\infty (\Bbb C P^\infty)$, again
denoted $1+k \psi ^{-k}$.

\vskip .2in
\proclaim{Theorem 2.2}
There are infinite loop maps
$$
\alpha : \Bbb Z\times B\Gamma ^+_\infty \lra \Omega ^\infty
\Bbb CP^\infty_{-1} , \quad
\mu_p: \Omega ^\infty S^\infty (\Bbb CP^\infty_+)
 \lra (\Bbb Z \times B\Gamma ^+_\infty) ^\wedge _p
$$
such that the composition
$\omega \circ \alpha \circ \mu_p$  and the self map
$$
\left( \matrix
1+ k \psi ^{-k}& 0  \\
0 & -2 \endmatrix \right)  :
\Omega ^\infty S^\infty \Bbb CP^\infty \times \Omega ^\infty S^\infty \lra
\Omega ^\infty S^\infty \Bbb CP^\infty \times \Omega ^\infty S^\infty
$$
become homotopic after $p$-adic completion.
\endproclaim

This   is an improvement on the main theorem of [MT] where the
map in the lower left corner had been left undetermined.
For our calculations in Section 4.3
we need this map however to be zero. 
A proof of Theorem 2.2  is given
in the final Section 5.

\vskip .2in
\noindent
{\bf Remark 2.3.}
The reader is referred to [BK] for the notion of $p$-adic completion
(also called $\Bbb F_p$-completion). For connected, compact CW-complexes $X$
and infinite
loop spaces $\Omega ^\infty E$ of finite type one has
$$
[X, (\Omega ^\infty E)^\wedge _p] = [X, \Omega ^\infty E]\otimes \Bbb Z_{p},
\quad
H^* ((\Omega ^\infty E ) ^\wedge _p; \Bbb Z) = H^* (\Omega ^\infty E;
\Bbb Z_{p}).
$$
Furthermore, note that
the homotopy class of the map $\alpha$ in Theorem 2.2 is uniquely determined
by its composition with
$\Bbb Z \times B\Gamma _\infty \to \Bbb Z\times B\Gamma ^+_\infty$.
Indeed since $\Omega ^\infty \Bbb CP^\infty _{-1}$ is an infinite loop space
the induced map
$$
[\Bbb Z \times B\Gamma^+_\infty , \,  \Omega ^\infty \Bbb CP^\infty _{-1}]
\lra [\Bbb Z \times B\Gamma_\infty , \,  \Omega ^\infty \Bbb CP^\infty _{-1}]
$$
is an isomorphism. This is a standard property of
the plus construction, cf. [B].

\vskip .1in
Before we give a detailed description of $\alpha$ in the next section, we state here the third result.

\vskip .2in
\proclaim{Theorem 2.4. ([MW])} The map $\alpha $ is a homotopy
equivalence.
\endproclaim

\vskip .4in
\subheading{3. Characteristic classes of surface bundles}

\vskip .2in
\subheading{3.1. Universal surface bundles}
The methods used in this and the surrounding papers do
not use the mapping class groups directly but rather the topological groups
of orientation preserving diffeomorphisms of surfaces.
We briefly review the correspondence.

Let $F_{g,b}$ be a connected surface of genus $g$
with $b$ boundary circles. We write
$\diff (F_{g,b}; \partial)$ for the topological
group of orientation preserving diffeomorphisms
that keep (a neighborhood of) the boundary pointwise fixed.
For $g\geq2$, results from [EE] and [ES] yield
$$
B\Gamma _{g,b} \simeq B\diff (F_{g,b}; \partial)
$$
so that $B\Gamma _{g,b}$ classifies diffeomorphism classes of
smooth fibre bundles
$\pi: E \to X$ with fibre $F_{g,b}$ and standard boundary
behavior:
$$
\partial E = X \times \coprod _{1} ^b S^1, \quad \pi | \partial E = \text{proj}_
X.
$$
Similarly,
$$
B\Gamma ^s_{g,b} \simeq B\diff (F_{g,b} ; \partial \coprod \{ x_1, \dots , x_s\}
)
$$
where $x_1, \dots, x_s$ are distinct interior points of
$F_{g,b}$. Take $s=1$. Since $\diff (F_{g,b}; \partial)$ acts transitively
on the interior of $F_{g,b}$,
$$
E (F_{g,b}) := E\diff
(F_{g,b}; \partial) \times _{\diff (F_{g,b}; \partial)}F_{g,b}
\simeq B\diff(F_{g,b}; \partial \coprod \{x\}) \simeq B\Gamma ^1_{g,b}.
$$
The forgetful map $\pi: B\Gamma ^1_{g,b} \to B\Gamma _{g,b}$ corresponds
to the universal smooth $F_{g,b}$ bundle
$$
F_{g,b} \lra E (F_{g,b})
\lra B\diff (F_{g,b}; \partial).
\tag 3.1
$$
The central extension (1.2) is classified by \lq \lq the differential at $x$",
$$
\diff (F_{g,b}; \partial \coprod \{x\}) \lra \text{GL} ^+ (T_x F_{g,b}) \simeq
SO(2).
$$
Hence the circle bundle induced from (1.2) by applying the classifying space
functor corresponds to the circle bundle of the vertical tangent
bundle associated with (3.1).

\vskip .3in
\subheading{3.2. The map $\alpha$ and the kappa classes}
The map of infinite loop spaces
$\alpha : \Bbb Z \times B\Gamma _\infty ^+\to \Omega ^\infty \Bbb C P^
\infty _{-1} $ constructed in section 2 of [MT] restricts to a map
$\alpha _{g,2}: B\Gamma _{g,2} \to
\Omega ^\infty _g \Bbb CP^\infty _{-1}$ that is homotopic to the
composition
$$
\alpha _{g,2} : B \Gamma _{g,2} \lra B\Gamma _{g+1} \overset
{\alpha _{g+1}} \to \lra \Omega ^\infty _g \Bbb C P^\infty _{-1},
\tag 3.2
$$
where the left hand map is induced from gluing the two
parametrized boundary circles together.
These maps are up to homotopy compatible with the monoidal structure
on $\bigcoprod B\Gamma _{g,2}$.

We next recall a description of $\alpha _{g+1}$ which
is well-suited for identifying the
kappa classes.
Let $\pi : E \to X$ be a smooth surface bundle with closed fiber $F$. Thus $E=
P \times _{\diff (F)}
F$ where $P$ is a principal $\diff (F) $ bundle over $X$. We
do not assume that $X$ is smooth or finite dimensional, only that $X$ is
paracompact (or a CW-complex).

We denote by $\emb (F, \Bbb R^n)$ the space of smooth embeddings in the
$C^\infty$-topology, and let $\Bbb R^\infty$ and
$\emb (F, \Bbb R^\infty)$ be the colimits of $\Bbb R^n$ and
$\emb (F, \Bbb R^n)$, respectively. We shall consider fiberwise embeddings
$\iota : E \to X \times \Bbb R^\infty$, that is, fiberwise maps such that each
$\iota _x: E_x \to \{ x\} \times \Bbb R^\infty$ is an embedding and such that the
adjoint $\diff (F)$-equivariant 
map $P \to \emb (F, \Bbb R^\infty )$ is continuous. Such
an $\iota $ is equivalent to a section of $P\times _{\diff (F)} \emb (F, \Bbb
R^\infty)$. Note that $\emb (F, \Bbb R^\infty)$ is contractible so that 
such a section always exists.

An embedding $\iota _x: F \hookrightarrow \Bbb R^{n+2}$ extends to a map
from the normal bundle $N^n \iota _x= \{ (p,v) | v \bot T_pF\}$ into
$\Bbb R^{n+2}$ by sending $(p,v)$ to $p+ v$. (Here we have identified $F$ with
its image under $\iota _x$.) We call the embedding
$\iota _x$ {\it fat} if this map restricts to an embedding of the unit disk bundle
$D(N^n \iota _x)$. The subspace of fat embeddings
$\emb ^f (F, \Bbb R^\infty) \subset \emb (F, \Bbb R^\infty)$ is
contractible, since the inclusion is a homotopy equivalence by the tubular
neighborhood theorem and since $\emb (F, \Bbb R^\infty)$ is
contractible by Whitney's embedding theorem.  A fibrewise fat embedding
$\iota : E\to X\times \Bbb R ^\infty $ is then a section of the
fibre bundle
$P\times _{\diff (F)} \emb ^f (F, \Bbb R^\infty)$.

Suppose first that $\iota : E \to X \times \Bbb R ^{n+2}$ is a fibrewise fat
embedding of codimension $n$. The Pontryagin-Thom construction associates
a \lq \lq collapse " map onto the Thom space of the fibrewise normal
bundle,
$$
c_{\pi, \iota } : X _+\wedge S^{n+2} \lra D(N_\pi ^n \iota )/ S(N^n_\pi \iota )
= \Th (N^n _\pi \iota ).
$$
We are particularly interested in its adjoint map $X\to \Omega ^{n+2} \Th
(N_\pi ^n \iota )$.

Let $G(2,n)$ be the Grassmann manifold of
oriented 2-dimensional subspaces of $\Bbb R^{n+2}$, and let $U_n$ and $U_n^\bot$
be the two  complementary universal 
bundles over it of dimension 2 and $n$, respectively.
The fat embedding $\iota $ induces bundle maps
$$
 T_\pi E \lra U_n, \quad N^n_\pi \iota  \lra U^\bot _n
$$
and a commutative diagram
$$
\CD
X_+ \wedge S^{n+2} @> c_{\pi,\iota }>>       \Th (N^n_\pi \iota )
        @> s>> \Th (T_\pi E \oplus N^n_\pi \iota )   \\
@|      @VVV    @VVV	\\
X_+ \wedge S^{n+2} @>c_n >> \Th (U_n^\bot)         @> s>>
        \Th (U_n\oplus U_n ^\bot).
\endCD
\tag 3.3
$$
In the general case of a fiberwise fat embedding $\iota : E \to X\times
\Bbb R ^\infty$, the base space $X$ is the colimit  of the
subspaces
$$
X_n:= \{ x \in X | \, \iota _x (E_x) \subset \{x \} \times \Bbb R^{n+2}\},
$$
and the diagram
$$
\CD
(X_n)_+ \wedge S^{n+2}  @>>>    \Th (U_n^\bot)  \\
@VVV    @VVV    \\
(X_{n+1})_+ \wedge S^{n+3}      @>>> \Th (U_{n+1}^\bot)
\endCD
$$
is commutative since $U^\bot_{n+1} |_{G (2,n)} = U^\bot _n$. Taking
adjoints we get
$$
\alpha _{\pi, \iota }: X \lra \colim _{n\to \infty} \Omega ^{n+2} \Th (U_n^\bot).
$$
Since $\emb ^f(F, \Bbb R^\infty) $ is contractible, all sections of
$P\times _{\diff (F) } \emb^f (F, \Bbb R^\infty)$ are homotopic, and
consequently the homotopy class $[\alpha _{\pi, \iota }]$ is independent of the
choice of $\iota $.
We will therefore from now on suppress the subscript $\iota $.

Realification gives a $(2n-1)$ connected map from $\Bbb CP^n$ into the oriented
Grassmannian $G(2,2n)$ covered by a bundle map $L^\bot _n \to U^\bot _{2n}$.
Thus $G(2, \infty) \simeq \Bbb CP^\infty$ and
$$
\Omega ^\infty \Bbb CP^\infty _{-1}  = \colim \, \Omega ^{2n+2}
\Th (L_n^\bot) \overset \simeq \to \lra \colim \, 
\Omega ^{2n+2} \Th (U_{2n}^\perp)
$$
is  a homotopy equivalence. Altogether we have a well-defined homotopy class
$$
\alpha _\pi : X \lra \Omega ^\infty \Bbb CP^\infty_{-1}.
$$
For $X= B\diff (F_{g+1}) \simeq B\Gamma _{g+1} $ this is the map
$\alpha _{g+1}$ of (3.2).

\vskip .1in
Let us check that the image of $\alpha _{g+1}$, and hence the image of
$\alpha _{g,2}$, lie in the $g$-component of $\Omega ^\infty
\Bbb CP^\infty _{-1}$, or equivalently that the composition
$$
\text{proj} \circ \omega \circ \alpha_{g+1}:
B\Gamma _{g+1} \lra \Omega ^\infty S^\infty (\Bbb CP^\infty _+)
\lra \Omega ^\infty S^\infty
$$
lands in the $-2g$ component (with identification of components
chosen at the end of Section 2.1). Consider (3.3) with $X$
a single point and $E=F_{g+1}$. The bottom row  in (3.3) is thus 
$ s\circ c_n: S^{n+2} \to G(2, 2n) _+ \wedge S^{n+2}$ 
and we need to compute the degree of the composition of this map
with the projection onto $S^{n+2}$. This degree is
given by the 
evaluation of the pullback 
of the generator of $H^{n+2} ( S^{n+2})$ on the fundamental class
$[S^{n+2}]$. Under the projection the fundamental class is pulled back to 
the Thom class of the trivial bundle  $U_n \oplus U_n^\perp$.
Writing  $\lambda_U$ for  the Thom  class of the vector bundle $U$, 
the degree  is thus  given by
$$
\aligned
< c^*_\pi s^* (\lambda _{U_n}\dot \lambda_{U_n^\perp}), [S^{2n+2}]> &=<
c^*_\pi (e(TF_{g+1})
\dot \lambda _{U_n^\perp}), [S^{2n+2}]>   \\
&=< e(TF_{g+1}) , [F_{g+1}] > = -2g,
\endaligned
$$
as claimed.

\vskip .1in
Next we compute the maps
$\alpha _{g,2}$ and $\alpha_{g+1,2}$ under the map $B\Gamma _{g,2} \to
B\Gamma _{g+1,2}$ induced from gluing a torus with two boundary
circles $F_{1,2} $ to $F_{g,2}$. Considering $F_{1,2}$ as a
fibre bundle over a point the construction above
gives an element $[1] \in \Omega ^\infty _1\Bbb C P^\infty _{-1}$.
Loop sum with $[1]$ in $\Omega ^\infty
\Bbb C P^\infty _{-1}$ translates the $g$ component into the $(g+1)$
component
and 
$$
\CD
B\Gamma _{g,2} @> \alpha _{g,2} >>
        \Omega ^\infty _g \Bbb CP^\infty _{-1}  \\
@VVV    @V *[1]VV       \\
B\Gamma _{g+1,2}        @>\alpha_{g+1,2}>>
        \Omega ^\infty _{g+1} \Bbb CP^\infty _{-1}
\endCD
\tag 3.4
$$
is homotopy commutative. 
To see this observe that 
the left vertical map is multiplication by the basepoint of $B\Gamma _{1,2}$
 in the monoid
$ \bigcoprod _{g\geq 0}
B\Gamma _{g,2}$ and the right vertical map is multiplication by its image in
 $\Omega ^\infty _1 \Bbb CP^\infty _{-1}$. Homotopy commutativity of 
diagram (3.4) now follows because  the $\alpha _{g,2}$ induce a map of
monoids up to homotopy for: $\alpha $ is a map of infinite loop spaces, 
$\alpha $ restricts to $\alpha_{g,2}$ on $B\Gamma _{g,2}$, and the infinite
loop space structure on $\Bbb Z\times B\Gamma ^+_\infty$ is compatible with
the monoidal structure.
(Alternatively,  homotopy commutativity of (3.4) follows from a calculations of
pretransfers, cf. [G2].)

\vskip .1in
Let $\tilde \alpha $ denote the restriction of $\alpha$ in Theorem 2.2 to the
zero component. Restricted to $B\Gamma _{g,2}$ it is homotopic to 
$\tilde \alpha _{g,2} = (*[-g])\circ\alpha_{g,2}.$
We can now relate the kappa classes to spectrum cohomology.
Consider
$$
B\Gamma _{g,2} \overset {\tilde \alpha _{g,2}} \to \lra
\Omega ^\infty _0 \Bbb CP^\infty _{-1}
\overset {\omega } \to \lra
\Omega ^\infty _0 S^\infty (\Bbb CP^\infty _+)
$$
and recall the cohomology suspension from Section 2.1
$$
\sigma ^*:  H^{2i} (\Bbb CP^\infty; \Bbb Z) \simeq H^{2i} (S^\infty
\Bbb CP^\infty _+) \lra H^{2i} (\Omega ^\infty _0 S^\infty (\Bbb CP^\infty_+)).
$$

\vskip .2in
\proclaim{Theorem 3.1}
The Miller-Morita-Mumford class $\kappa_i$ is equal
to $( \omega \circ \tilde \alpha )^* (\sigma ^*e^i)$  where
$e\in H^2 (\Bbb CP^\infty ; \Bbb Z)$ is the Euler class of the canonical
line bundle.
\endproclaim

\vskip .2in
\noindent
{\it Proof:}
Let $\pi : E \to X$ be a smooth fibre bundle with fibre
$F_{g+1}$, classified by $f_\pi : X \to B\Gamma _{g+1}$.
By definition
$$
f^*_\pi(\kappa_i) = \pi_!(e(T_\pi E)^{i+1}) \in H^{2i}(X; \Bbb Z)
$$
where $\pi_!$ is the composition of the Thom isomorphism and the
Pontrjagin-Thom collapse map
$$
\tilde H^{2i+2}(E; \Bbb Z) \overset \simeq \to \lra \tilde
H^{2i+2n+2}(\Th (N_\pi^n \iota ); \Bbb Z) \overset {c_{\pi}^*} \to \lra
H^{2i+2n+2} (S^{2n+2}\wedge X_+; \Bbb Z)
$$
followed by the $(2n+2)$-nd desuspension; here the notation  is as in
(3.3). 

Let $x \in H^k (\Bbb CP^N; \Bbb Z) \simeq H^k (G(2, 2N); \Bbb Z)$
for $N \gg k$. The 
$2n+2$-fold suspension $\Sigma ^{2n+2} (x) \in H^{k+2n+2} ( S^{2n+2}\wedge
\Bbb CP^N, \Bbb Z)$
is $x$ times the Thom class of the
trivial $2n+2$ real bundle $L_n^\perb \oplus L_n$. Thus 
$$
s^* ( \Sigma ^{2n+2} (x) ) = s^* ( \lambda _{L_n ^\perb} \dot \lambda _{L_n}
\dot x ) 
= \lambda_{L_n ^\perb} \dot e(L_n)\dot x.
$$
To interpret these formulas, recall that the cohomology of the Thom space
of a vector bundle is a module over the cohomology of the base space, and a
map of bundles  (such as $s$) induces a  map of modules. 
Furthermore, as elements in the cohomology of the Thom space
of the trivial bundle, $\lambda _{L_n ^\perb}$ and $\lambda _{L_n}$ are 
pulled back along the bundle projections $\pi _1$ and $\pi _2$ of $L_n^\perb
\oplus L_n$ onto the first and second summand. 
Hence, as $\pi_1 \circ s$ is the identity, we have $s^* ( \lambda _{L_n ^\perb})
=\lambda _{L_n^\perb}$, and as $\pi_2 \circ s$ factors through the base space,
$s^* (\lambda _{L_n}) = 1\dot e (L_n)$.
Finally, take $x=e^i$ in the  above formula and use the commutativity
of (3.3) to complete the proof.
\qed

\vskip .3in
\subheading{3.3. One part of Akita's conjecture}
For our next theorem we need the relation between Steenrod operations and
characteristic classes of vector bundles.
Recall
the $i$-th Steenrod operation:
$$
\aligned
&P^i : H^k(X; \Bbb F_p) \lra H^{k+2i(p-1)}(X; \Bbb F_p),\quad p \text{ odd},
\\
& Sq^i: H^k (X, \Bbb F_2) \lra H^{k+i} (X, \Bbb F_2).
\endaligned
$$
Let $E$ be an oriented vector bundle over $X$ and $\lambda _U$
its cohomology  Thom class. One defines $v_i (E) \in H^{2i(p-1)}
(X; \Bbb F_p)$, respectively $v_i (E) \in H^i (X; \Bbb F_2)$ by
$$
P^i (\lambda _E) = v_i (E) \lambda _E, \quad
Sq^i (\lambda _E) = v_i (E) \lambda _E.
$$
For $p=2$, these are  the Stiefel-Whitney classes, and for $p$
odd they were first defined by Wu, cf. [MSt].
For an oriented 2-plane bundle (or complex line bundle) $L$,
$$
\aligned 
v_1 (L) = e(L)^{p-1} \quad \quad \quad    \quad &\text{ for $p$ odd}, \\
v_2 (L) = e(L) \quad \text{and} \quad v_1 (L) = 0 \, \quad &\text{for } p=2.
\endaligned 
$$
Moreover, the total class
$$
v (E) = 1 + v_1 (E) + v_2 (E) + \dots \in H^* (X; \Bbb F_p)
$$
takes direct sums of oriented vector bundles into (graded) products.

\vskip .2in
\proclaim{Theorem 3.2}
The modulo $p$ reduction of $\kappa_i \in H^{2i} (B\Gamma _\infty;
\Bbb F_p)$ is zero when $i+1 \equiv 0 \, (\text{mod } p-1)$.
\endproclaim

\vskip .2in
\noindent
{\it Proof:}
\footnote{
We thank John Rognes for this proof; it replaces a more
cumbersome earlier argument.}
Let $H(\Bbb Z, k)$
denote the Eilenberg-MacLane space with non-trivial homotopy
$\Bbb Z$ in dimension $k$. The Thom class $\lambda_n
 = \lambda_{L_n ^\perb}$ is represented by a map from
$\Th (L_n^\perb)$ to $H(\Bbb Z, 2n)$. We let $Y_{2n+2}$ be its homotopy fibre,
so that there is a fibration sequence
$$
Y_{2n+2} \overset {j_{2n+2}} \to \lra \Th (L_n^\perb) \overset \lambda_n \to
\lra H(\Bbb Z, 2n).
$$
The spaces $Y_{2n+2}$ are the $(2n+2)$-nd terms of a spectrum $Y$ and the
$j_{2n+2}$ define a map $j: Y \to \Bbb CP^\infty _{-1}$ of spectra.
Since $\Omega ^{2n+2} H(\Bbb Z, 2n)$ has vanishing homotopy groups
$$
\Omega ^\infty Y \overset {\Omega ^\infty j} \to \lra
\Omega ^\infty \Bbb CP^\infty _{-1}
\tag 3.5
$$
is a (weak) homotopy equivalence.
Thus  the bottom vertical map in the commutative diagram
$$
\CD
H^{2i} ( \Bbb CP^\infty _{-1}; \Bbb F_p) @ > j^* >> H^{2i} (Y, \Bbb F_p)
\\
@V \sigma ^* VV @V \sigma ^*VV  \\
H^{2i} (\Omega ^\infty _0 \Bbb C P^\infty _{-1} ; \Bbb F_p)
@> ( \Omega ^\infty j )^* >> H^{2i} (\Omega ^\infty Y; \Bbb F_p)
\endCD
$$
is an isomorphism.
The proof of Theorem 3.1 shows that
$$
\tilde \alpha ^* \sigma ^* (e^{i+1} \lambda _{L^\perb}) = \kappa_i,
$$
so it suffices to prove that $j^*$ vanishes when $i+1 \equiv 0
\, (\text{mod}\, \,  p-1)$.
Equivalently, we must show that $e^{i+1} \lambda _n$ is in the image of
$$
\lambda ^* _n: H^* (H(\Bbb Z, 2n); \Bbb F_p) \lra H^* (\Th (L_n ^\perb); \Bbb
F_p)
$$
in the stated dimensions. This is implied by
$$
v(L_n ^\perb) = v(L_n )^{-1} = \left\{
        \aligned
         (1 + e ^{p-1}) ^{-1} , \quad \quad  & p>2     \\
         (1+ e)^{-1}, \quad \quad & p=2
        \endaligned
	\right.
$$
Since $P^i (\lambda _n) \in \text{image } (\lambda_n^*)$ the result follows.
\qed

\vskip .2in

\subheading{ Remark 3.3}
The relation $P^i (\lambda _{L^\perb}) = \kappa _{i(p-1) -1} \lambda _{L^\perb}$
used above is further exploited in [G2] to define secondary classes $\mu_i$
with $p \mu_i = \kappa _{i(p-1) -1}$ in cohomology with $\Bbb Z / p^2\Bbb Z$
coefficients.

\vskip .1in
Theorem 3.2 proves half of Akita's conjecture mentioned in the introduction.
The other half is implied by Theorem 1.1 which is proved below.

\vskip .4in
\subheading{ 4.  Proof of Theorem 1.1}

\vskip .2in
\subheading{4.1. Segal's splitting}
The inclusion of $\Bbb CP^\infty $ in $BU$ that represents the reduced canonical
line bundle (of virtual dimension zero) extends to a map
of infinite loop spaces
$$
l: \Omega ^\infty S^\infty (\Bbb CP^\infty ) \lra 
 BU
$$
by Bott periodicity and
the universal property of the free infinite loop space functor 
$\Omega ^\infty S^\infty$.
Graeme Segal [S] proved that this map has a splitting. In [C] Michael
Crabb gives a construction of this map as 
the $S^1$-equivariant $J$-homomorphism.

\vskip .2in
\proclaim {Theorem 4.1 ([S], [C])} 
The map $l$ has a left inverse up to homotopy. In the resulting decomposition
$$
\Omega ^\infty S^\infty (\Bbb C P^\infty) \simeq BU \times Fib(l)
$$
the homotopy fiber $Fib(l)$ has vanishing rational cohomology.
\endproclaim

%
%

\vskip .2in
In particular, this gives an identification of Hopf algebras
$$
\FH^*(\Omega ^\infty S^\infty (\Bbb C P^\infty)) \simeq H^* (BU; \Bbb Z)
\tag 4.1
$$
induced by $l^*$.
The graded module of primitive elements of the right hand side
is a copy of $\Bbb Z$
in each even degree generated by the integral Chern character class
$s_i = i!ch_i$. 
Since $l\circ i$ represents the (reduced) line bundle and $ch_i (L) = \frac
1 {i!} e^i $, (4.1) implies that 
$$
i^*: \FH^* (\Omega ^\infty S^\infty (\Bbb CP^\infty)) \lra
H^*(\Bbb CP^\infty; \Bbb Z)
$$
sends the primitive generator
$s_i$ to $e^i$. Note from (2.2) that we also have that
$$
H^* (\Bbb CP^\infty ; \Bbb Z)= H^* (S^\infty \Bbb CP^\infty ; \Bbb Z)
\overset {\sigma ^*} \to \lra P(\FH^* (\Omega ^\infty S^\infty 
(\Bbb CP^\infty)))
$$
maps $e^i$ to $s_i$.
Therefore, Theorem 3.1 translates into
$$
\kappa _i = (l\circ \omega \circ \alpha )^* (s_i).
\tag 4.2
$$

\vskip .1in
The map $\omega$ in (2.3) induces an isomorphism on rational cohomology because the cohomology of $\Omega ^\infty S^{\infty -1}$  is well-known to be all
torsion.
The isomorphisms
$$
H^*(BU; \Bbb Q) \overset {l^*} \to \lra H^*(\Omega ^\infty_0 S^\infty
(\Bbb CP^\infty _+); \Bbb Q) \overset {\omega ^*} \to \lra
H^* (\Omega ^\infty \Bbb CP^\infty _{-1}; \Bbb Q)
$$
together with
Theorem 2.4 thus
give an isomorphism
$$
H^*(B\Gamma_\infty ; \Bbb Q) \simeq H^*(BU; \Bbb Q)
$$
of Hopf algebras.  Hence  $P(\FH^* (B\Gamma _\infty))$ is a copy of $\Bbb Z$
in each even degree. We choose a generator
$$
\tau_i \in P(\FH^{2i} (B\Gamma_\infty))
$$
such that $\kappa_i = D_i \tau_i$ for some positive number $D_i$.
The next section gives close upper and lower bounds for $D_i$.

\vskip .2in
\subheading{4.2. A  lower  and an upper bound}
Recall the definition of the symplectic characteristic classes for 
surface bundles.
The action of $\Gamma _{g,1} $ on $H^1 (F_{g, 1}; \Bbb Z) = \Bbb Z^{2g} =
H^1 (F_{g}; \Bbb Z)$
induces the standard symplectic representation, and hence a representation of 
$\Gamma_{g,2}$ via the map $\Gamma_{g,2} \to \Gamma_{g, 1}$ induced by
gluing a  disk to one of the boundary components.
We may let $g\to \infty$ and obtain
$$
B\Gamma_\infty \lra B\text{Sp}(\Bbb Z).
$$
This  map can be composed with the map into $B\text{Sp}(\Bbb R)
 \simeq BU$ so that we have a map
$$
\eta: B\Gamma _\infty \lra BU.
$$

\vskip .2in
\proclaim{Theorem 4.2 ([Mo], [Mu])}
In $H^*(B\Gamma_\infty ; \Bbb Q )$ one has the relation
$$
\eta^* (s_{2i-1}) = (-1)^i (\frac {B_i}{2i})\, \kappa_{2i-1}.
$$
\endproclaim

We are now in the position to prove

\vskip .2in
\proclaim{Theorem 4.3} For all $i\geq 1$, 
$D_{2i} =2 $  and 
$$
D_{2i -1} = \den ( \frac { B_i} {2i}) 
\quad \quad \text { or } \quad \quad  D_{2i-1} 
=  2 \,  \den (\frac{B_i} {2i}).
$$
\endproclaim

\vskip .2in
\demo{Proof}
As $s_{2i-1}$ and $\kappa_{2i-1}$ are integral classes, Theorem 4.2
implies immediately
that modulo torsion $\den (B_i/2i)$ divides $\kappa_{2i-1}$.
By  Theorem 3.2 (for $p=2$) we also know that $2$ divides $\kappa _{2i}$.
This establishes the lower bounds for all $i \geq 1$.
The upper bounds
are a consequence of Theorem 2.2 as we explain now.

\vskip .1in
As in Section 2.2, let $k$ be   a positive 
integer such that $-k$ generates $(\Bbb Z/p^2) ^\times$ for odd $p$ and 
let $k=5$ when $p=2$. Then by Theorem 2.2
there is a factorization
$$
1+k \psi^{-k} : \Bbb CP^\infty \overset \mu_p \to \lra
(B\Gamma ^+_\infty)^\wedge _p \overset {\omega \circ \alpha} \to \lra
\Omega ^\infty _0 S^\infty (\Bbb CP^\infty _+) ^\wedge _p
\overset \text{proj} \to \lra \Omega ^\infty S^\infty (\Bbb CP^\infty) ^\wedge
_p.
$$
Since $\Omega ^\infty S^\infty (\Bbb CP^\infty)$ is of finite type
$$
\FH^*(\Omega ^\infty S^\infty (\Bbb CP^\infty) ^\wedge _p) = 
\FH^* (\Omega ^\infty S^\infty (\Bbb CP^\infty ))\otimes \Bbb Z_{p}.
$$
It follows from (4.2)
that $\mu^*_p (\kappa_i)$ is the image of the primitive
generator $s_i$ under
$$
(1+k \psi^{-k} )^*: P\FH^{2i} (\Omega ^\infty S^\infty (\Bbb CP^\infty) )\otimes
\Bbb Z_p  \lra H^{2i} (\Bbb CP^\infty ; \Bbb Z_p).
$$
Both groups are copies of $\Bbb Z_p$ and $(1+k\psi^{-k})^*$ in  dimension $2j$
is multiplication
by
$1+k(-k)^j = 1- (-k) ^{j+1}$.
Therefore, $\nu_p (D_j) \leq \nu _p ( 1-(-k)^{j+1})$.

\vskip .1in
We have the following well-known table of
$p$-adic valuations (see for example lemma 2.12 in [A2]):
$$
{\nu _p (1 - (-k)^s) = \left\{ \aligned
        1+\nu _p (s) \quad \quad &\text { if } s \equiv 0 \quad \mod (p-1), 
		\, \, p \text{ odd }\\
        0 \quad \quad &\text{ if } s   \nequiv 0 \quad \mod (p-1),
		\, \, p \text { odd }
        \endaligned \right. } 
$$
$$
{\nu _2 (1 - (-k)^s) = \left\{ \aligned
        2+\nu _p (s) \quad \quad &\text { if } s \equiv 0 \quad \mod (2)\\
        1 \quad \quad &\text{ if } s  \nequiv 0 \quad \mod (2),
        \endaligned \right. }
%
$$
When $j=2i$,  for $p$ odd $\nu_p (1-(-k)^{j+1}) $ is  zero  and for $p=2$
it 
is equal to 1. This gives $D_{2i} \leq 2$.
When $j=2i-1$,  for $p$ odd   we have
$\nu _p (1- (-k) ^{j+1}) = 1 + \nu _p (2i)$ if $2i$ divides $p-1$ and zero
otherwise. This  
is precisely $\nu _p (\den (B_i / 2i))$ (compare (1.3)).
If $p=2$ however, $\nu _2 (1- (-k)^{j+1}) = 2 + \nu _2 (2i) $ which is
one more than $\nu_2 (\den (B_i /2i))$. 
This gives $D_{2i-1} \leq 2 \, \den (B_i /2i)$.
\qed 
\enddemo

\vskip .3in
\subheading{4.3. The final factor of 2}
Theorem 4.3 leaves us with an indeterminancy of a factor of 2 in the 
odd case, and Theorem 1.1
will follow immediately from Theorem 4.5 below.
The proof requires several extra results: 
the improvement of the main theorem from [MT] as stated in Theorem
2.2 (for $p=2$), 
Theorem 2.4, as well as  part of theorem 1.3 from [G1] which we state
as

\vskip .2in
\proclaim{Theorem 4.4}
$
H_* (\Omega ^\infty _0 \Bbb CP^\infty _{-1}; \Bbb F_2) \overset {\omega _*}
\to \lra H_* (\Omega ^\infty _0 S^\infty (\Bbb CP^\infty _+); \Bbb F_2)
$
is injective.
\endproclaim

\vskip .2in
\proclaim{Theorem 4.5}
$
\nu _2 (D_{2i-1}) = 1+\nu_2 (2i).
$
\endproclaim

\vskip .1in
\noindent
{\it Proof:}
As in the proof of
Theorem 4.3  we have
$$
\aligned
1- (-5) ^{j+1} : P(\FH^{2j} (\Omega ^\infty S^\infty \Bbb CP^\infty _+)))
\otimes
\Bbb
Z_2 \overset {(\omega \circ \alpha )^* } \to \lra 
&P(\FH^{2j} (B\Gamma ^+_\infty))
\otimes \Bbb Z_2        \\
 \overset {\mu_ 2^*}
\to  \lra  &H^{2j} (\Bbb CP^\infty; \Bbb Z_2).
\endaligned
$$
All groups are copies of $\Bbb Z_2$, $\kappa_j =D_j\tau_j$ with $\tau_j$
a generator and $\kappa_j = (\omega \circ \alpha )^* (s_j)$ where $s_j$
is the generator of the left term.
Suppose that $\nu_2 (D_j) =2 +\nu _2(j+1)$ for some $j = 2i-1$.
Then $\mu^*_2(\tau_j)$ would be a generator and dually
$$
{\mu_2}_* : H_{2j} (\Bbb CP^\infty ; \Bbb F_2) \lra
H_{2j}^{free} (B\Gamma ^+_\infty)
\otimes \Bbb F_2
$$
and hence
$$
{\mu_2}_* : H_{2j} (\Bbb CP^\infty ; \Bbb F_2) \lra H_{2j} (B\Gamma ^+_\infty)
\otimes \Bbb F_2
$$
would  be non-zero. Now apply Theorem 2.4 and Theorem 4.4 
 to conclude
that
$$
H_{4i-2}(\Bbb CP^\infty; \Bbb F_2) \overset (\omega \circ \alpha \circ
\mu_2 )_* \to \lra H_{4i-2} (\Omega ^\infty _0 S ^\infty (\Bbb CP ^\infty _+);
\Bbb F_2)
$$
would be non-zero.
But this leads to a contradiction as we now argue.

\vskip .1in
Indeed, by Theorem 2.2 the above map $(\omega \circ \alpha \circ \mu_2)_*$ is 
$(1 + 5\psi ^{-5}, 0)_*$.
The self map $\psi^{-5}$ of $\Bbb
CP^\infty$ induces the identity on $H_* (\Bbb CP^\infty; \Bbb F_2)$. Hence
$1+5\psi^{-5}$  induces the same map on $\Bbb F_2$-homology as
six times the canonical inclusion map
$$
6i: \Bbb CP^\infty \lra \Omega ^\infty S^\infty (\Bbb CP^\infty).
$$
But $(6i)_*=0$ on $\Bbb F_2$-homology in degrees $4i -2 $ giving the desired
contradiction.
Indeed, $6 i$ is the composition of $2i$ and multiplication by $3$ in the
loop sum sense in $\Omega ^\infty S^\infty (\Bbb CP^\infty)$.  
$2i$ in turn  is the composition
of the diagonal map on $\Bbb CP^\infty$ composed with $i\times i$ and 
loop sum. The diagonal map
sends the generator $a_j$ of the $2j$-th homology group to 
$\Sigma _{s+t = j} a_s \otimes a_t$. Loop sum replaces tensor product 
by Pontryagin product which  here is commutative. 
Thus in $\Bbb F_2$-homology for $j$ odd, $2i$ is zero. 
\qed

\vskip .4in
\subheading {5. Proof of Theorem 2.2}

\vskip .2in
We offer two proofs of Theorem 2.2. The first one is in the spirit of [MT] 
and provides an
improvement on the geometric construction of $\mu_p$ given there.
The second
proof was suggested by the referee.  It is short, purely homotopy
theoretic, but it does not provide any insight in the map $\mu_p$ from
$\Omega^\infty S^\infty \Bbb CP^\infty_+$ to $B\Gamma_\infty^+$.  The map
$\mu_p$, constructed from explicit cyclic ramified covers of $\Bbb CP^1$,
provides an important link between the algebraic approach and the
topological approach to the mapping class group.

\vskip .2in
\subheading{5.1. A constructive proof}
In order to prove that the map  $\Omega ^\infty S^\infty \Bbb CP^\infty \to
\Omega ^\infty S^\infty$  induced by $\omega \circ \alpha \circ \mu_p$ is
zero after $p$-adic completion, we will have to review the construction
of $\mu_p$.  $\mu_p$ depends on  $p$ and the choice of $k$;
here $k$ can be any integer which is  a unit modulo $p^n$ for all $n
\geq 1$.  We take  the opportunity to give a variant of the proof
for the main theorem in [MT] from sections 3.1-3.3.  
For the fact that $\alpha $ is a map of infinite loop spaces we refer to 
section 2 of [MT]. 

\vskip .1in
We first describe Riemann surfaces $\Sigma$ with holomorphic actions of the 
$q$-th roots of unity $\mu_q \subset \Bbb C^\times$. 
This gives maps $B\mu_q \to B\diff (\Sigma)$. 

Consider a divisor $D= \Sigma  n_i p_i $ of $\Bbb CP^1$ with support 
$A = \{ p_0, p_1, \dots, p_k\}, n_i \in \Bbb N$ and  $n_0 + n_1 + \dots + n_k
\equiv 0 (\mod q)$. 
Assume for simplicity that gcd$(q,n_i)= 1$ for
$ i=0, 1,\dots , k$.
Let $\Sigma _D$ be the branched cover associated with the Galois extension
$$
\Bbb C (z) \hookrightarrow \Bbb C (z)[T]/ (F(T)), \quad \quad 
F(T) = T^q - \Pi_{ i=0} ^k (z- p_i) ^{n_i}
$$
(see e.g. [F], chap. 1.8). The Galois group is the group $\mu_q$. The surface 
$\Sigma _D$ has a holomorphic action of  $\mu_q$  with orbit space $\Bbb CP^1$.
The induced map $\pi : \Sigma _D \to \Bbb CP^1$ is holomorphic, branched over
$A$, and $\pi ^{-1} (p_i)$ is a single point 
for each $i$ (since we assumed gcd$(n_i, q) =1$).
Thus the $\mu_q$ action on $\Sigma _D$ is free outside $A$, and $A$ is fixed 
pointwise by all elements in $\mu_q$.

Let $\gamma_i $ be a small loop in $\Bbb CP^1$ around $p_i$. The fundamental group $\Bbb CP^1 \setminus A$ is the free group of rank $k$ generated by $\gamma_0,
\dots , \gamma _k$ with the single relation $\Pi \gamma _i=1$. The covering 
$\Sigma _d \setminus A \to \Bbb CP^1 \setminus A$ is classified by the map
from
$\pi _1 (\Bbb CP^1 \setminus A)$ to $\mu_q$ that sends $\gamma_i$ to
$e^{2\pi i n_i/ q}$. The complex tangent line $T_{p_i}\Sigma _D$ at $p_i$ is
a $\mu_q$ representation; $u \in \mu_q$ multiplies by $ u ^{\bar {n_i}}$ 
where $\bar n_i \in \Bbb Z/q$ is the multiplicative inverse of $n_i$.

If $D$ and $D'$ are two divisors and $q$ divides their difference $D-D'$ 
then there is a biholomorphic map between $\Sigma _D$ and $\Sigma _{D'}$ that 
is equivariant w.r.t. the $\mu_ q$ action. Thus it is only the class of $D$ in
$\tilde H_0 (A; \Bbb Z /q)$ that matters. 

Recall from Section 3.2 
that given a surface bundle $\pi : E \to X$ there is a diagram
$$
\CD
X @>>> \Omega ^\infty (\Th (N_\pi E)) 
	@>>> \Omega ^\infty S^\infty (E_+)	\\
@|	@VVV	@V T_\pi VV	\\
X @> \alpha _{p_i} >> \Omega ^\infty \Bbb CP^\infty _{-1} @>\omega >> 
\Omega ^\infty S^\infty  (\Bbb CP^\infty _+)
\endCD
\tag 5.1
$$
where $\Omega ^\infty ( \Th (N_\pi E)) = \colim \, \Omega ^{n+2} (N_\pi ^n
e_n)$ for suitable fat fibrewise embeddings $e_n : E
\to X \times \Bbb R^{n+2}$. 
The upper horizontal composition is the Becker-Gottlieb transfer map
$ t_E =t_\pi: X  \to \Omega ^\infty S^\infty (E_+)$.

We need the following properties of the Becker-Gottlieb transfer for smooth
manifold bundles with compact fiber and compact Lie structure group:

\roster
\item "(A1)" Let $f : E \to E'$ be a fiberwise homotopy equivalence. Then
$$
t_{E'}= \Omega ^\infty S^\infty (f_+) \circ t_E \in
[X, \Omega ^\infty S^\infty (E'_+)].
$$

\item "(A2)" Suppose 
$$
\CD
E_{12} 	@> j_1>>	E_1	\\
@V j_2 VV	@V i_1VV	\\
E_2 	@> i_2>>	E
\endCD
$$
is fiberwise homotopy coCartesian. In $[ X, \Omega ^\infty S^\infty(E_+)]$,
$$
t_E = \Omega ^\infty S^\infty ({i_1}_+) \circ t_{E_1} 
+ \Omega ^\infty S^\infty ({i_2}_+) \circ t_{E_2} -
\Omega ^\infty S^\infty ({i_{12}}_+) \circ t_{E_{12}}
$$
where 
$i_{12}= j_1\circ i_1= j_2\circ i_2$.

\item "(A3)" If the tangent bundle along fibers $T_\pi E$ admits a non-zero 
section, then $t_E \in [ X, \Omega ^\infty  S^\infty (E_+)]$ is the zero
element.

\endroster

The proof of (A1) and (A2) can be found in [LMS], p. 189-190 or in [BS].
Property (A3) is much simpler. It follows because 
$$
\Th (N_\pi E) \lra \Th( N_\pi E \oplus T_\pi E)
$$
is homotopic to the constant map at $\infty$ whenever $T_\pi 
E$ has an everywhere
non-zero section.

Let $\Sigma =\Sigma _D$ be the $\mu_q$-surface constructed above.
We shall study the transfer of the associated smooth surface bundle
$$
\pi : E\mu_q \times _{\mu_q} \Sigma \lra B\mu_q.
$$
To shorten notation we write 
$$
t_\Sigma : B\mu_q \lra \Omega ^\infty (E_+), \quad
 E= E\mu_q \times _{\mu_q} \Sigma
$$
for the associated transfer.

Let $D= \Sigma n_i p_i $ and $q \in \Bbb N$ satisfy 
$A= \text{ supp} (D) = \{ p_0, \dots, p_k\}, n_0 + \dots +n_k \equiv 0,$
gcd$(q,n_i)=1$ for $i=0, \dots, k$.
For each $i$, the inclusion of $E\mu_q \times _{\mu_q} \{p_i\} \subset
E$ induces a map 
$$
\hat p_i: B\mu_q \lra E \lra \Omega ^\infty S^\infty (E_+).
$$
The principal $\mu_q$ bundle $E\mu_q \to B\mu_q$ induces a transfer from
$B\mu_q$ to 
$\Omega ^\infty S^\infty (E{\mu_q}_+) $ $ \simeq \Omega ^\infty S^\infty$,
and hence 
$$
\hat t _q : B\mu_q \lra \Omega ^\infty S^\infty 
\lra \Omega ^\infty S^\infty (E_+)
$$
upon choosing a point of $E$.

\vskip .2in
\proclaim{Lemma 5.1}
The transfer $t_\Sigma $ is equal to $\Sigma \hat p_i + (1-k) \hat t_q$
in $[ B\mu_q, \Omega ^\infty S^\infty (E_+)]$.
\endproclaim

\demo{Proof}
We make a cell decomposition of $S^2 =\Bbb CP^1$ with two $0$-cells 
$\{ 0, \infty\}$, $k+1$ 1-cells $I_i$ and $k+1$ 2-cells $D_i$ such that $p_i
\in \text{ int }  D_i$.

\vskip .3in
\centerline{\hbox{\hskip 0.2cm \epsfysize=5cm\epsfbox{globe.eps}}}

\vskip .2in
\hfil {\it  Figure 5.2} \hfil

\vskip .2in
There are obvious coCartesian diagrams
$$
\CD
\coprod \partial I_i @>>> \coprod I_i @. \quad\quad \quad\quad 
	 \coprod \partial D_i @>>> \coprod D_i\\@VVV	@VVV
	\quad\quad	@VVV 	@VVV \\
\{0, \infty\} @>>> G @. \quad \quad  G @>>> S^2
\endCD
$$
where $G$ denotes the 1-skeleton. This cell structure lifts to a cell structure of $\Sigma $:
$$
\CD
\coprod \partial \tilde I_i @>>> \coprod \tilde I_i @. \quad \quad \quad \quad 
\coprod \partial \tilde  D_i @>>> \coprod \tilde D_i\\
@VVV    @VVV  \quad \quad  @VVV    @VVV \\
\{0, \infty\}^\sim  @>>> \tilde G @. \quad \quad \tilde G @>>> \Sigma
\endCD
\tag 5.3
$$
with $\tilde G = \pi ^{-1} (G)$ etc. We can apply the functor 
$E\mu_q \times _{\mu_q} (-)$ to (5.3) 
and use (A1-3) to evaluate $t_{\tilde G}$ and $t_\Sigma$. First by (A1) and (A2),
$$
t_{\tilde G} = (1+k) \hat t_q + 2 \hat t _q - 2(k+1) \hat t _q = (1-k) \hat t_q,$$
since $\tilde I_i = \mu _q \times I_i, \partial \tilde I _i = 
\mu _q \times \partial I_i$ and $ \{ 0, \infty \} ^\sim = \mu _q \times \{0
, \infty \}$.
Second, the inclusion of $p_i$ in $ \tilde D_i$ is a homotopy equivalence,
so $t_{\tilde D_i }= \hat p_i$. Moreover,
$ E\mu_q \times_ { \mu_q} \partial \tilde D_i= E\mu _q \times S^1$ has trivial 
vertical tangent bundle, so $t_{\partial \tilde D_i}$ is homotopically constant.
One more application of (A1-2) completes the proof.
\qed
\enddemo

\vskip .2in
The tangent representation $T_{p_i} \Sigma $ is given by multiplication with
$e ^{ 2\pi i \bar n_i/q}$, so an application of Lemma 5.1 gives

\vskip .2in
\proclaim{ Corollary 5.2}
The homotopy class of
$$
B\mu _q \lra \Omega ^\infty \Bbb CP^\infty _{-1} \lra \Omega ^\infty
 S^\infty (\Bbb CP^\infty _+)
$$
is $\Sigma \psi ^{\bar n_i} + (1-k) \hat t_q$. Here 
$\psi ^{\bar n_i}$ is the composition
$$
\psi ^{\bar n_i}: B\mu _q \lra \Bbb CP^\infty \lra \Omega ^\infty S^\infty
(\Bbb CP^\infty _+)
$$
with the left hand map induced from the group homomorphism 
$\mu_q \to S^1$ that sends $u$ to $u^{\bar n_i}$.
\qed
\endproclaim

\vskip .2in
We are now ready to complete the proof of Theorem 2.2. As in section 3.3 of 
[MT] we let $q=p^n$ be a prime power and consider the divisor
$$
D= p_0 +m p_1 + \dots +mp_k, \quad \quad  m\equiv -1/k (\text{mod } p^n).
$$
We use the notation 
$$
F(n)= \Sigma _D, \quad C_{p^n} = \mu _{p^n}, \quad \tau _n =\hat t _{p^n}
$$
and consider diagram (5.1) with
$$
X= B \diff (F(n)), \quad \quad  E = E \diff (F(n)) \times _{\diff (F(n))} F(n).
$$
Composing with the map $BC_{p^n} \to B\diff (F(n))$ induced by the $C_{p^n}$
action on $F(n)$, and using that $B\Gamma _{g(n)} \simeq B \diff (F(n))$ 
with $g(n) = 1/2 (p^n -1) (k-1)$
we get the diagram
$$
\CD
B\Gamma _g(n) @> \alpha _n>> \Omega ^\infty _{g(n) -1} \Bbb CP^\infty _{-1}
	@> T>>	\Omega ^\infty _0 \Bbb CP^\infty _{-1}	\\
@A\mu_n AA @ V \omega VV @V\omega VV	\\
BC_{p^n} @>>> \Omega ^\infty _{2- 2g(n)} S^\infty (\Bbb CP^\infty _+) 
	@ > T >> \Omega ^\infty _0 S^\infty (\Bbb CP^\infty _+).
\endCD
$$
The right-hand horizontal  
maps are translations of the indicated component into the zero component. 
The  lower horizontal composition is by Corollary 5.2 equal to 
$$
(1+k\psi^{-k} , \tilde \tau _n) \in [BC_{p^n}, \Omega ^\infty  S^\infty
 (\Bbb CP^\infty ) \times \Omega ^\infty _0 S^\infty]
$$
with $\tilde \tau _n = T \circ \tau _n$.

\vskip .2in
\proclaim{Lemma 5.3}
Let $i_{n-1}: BC_{p^n-1} \to BC_{p^n}$ be the map associated with $C_{p^{n-1}}
\subset C_{p^n}$. Then 
$$
[ \tau _n \circ i_{n-1}] = p [\tau _{n-1}] \in [BC_{p^{n-1}}, \Omega ^\infty
S^\infty].
$$
\endproclaim

\vskip .2in
\noindent
{\it Proof:}
Let $E$ be a contractible space with a free action of $C_{p^n}$ for all $n$,
e.g. the union of odd dimensional spheres $E= \bigcup _{m\geq 1} S^{2m-1}$.
Consider the diagram
$$
\CD
EC_{p^n} @> \pi_{n-1}>> BC_{p^{n-1}} @> i_{n-1} >> BC_{p^n}	\\
@AAA	@AAA	@A i_{n-1}AA	\\
\bigcoprod ^p_1 EC_{p^n} @>>> \bigcoprod ^p_1 BC_{p^{n-1}} @>>> BC_{p^{n-1}}
\endCD
$$
where $BC_{p^{n-1}}$ and $BC_{p^n}$ are the orbit spaces $E/C_{p^{n-1}}$
and $E/C_{p^n}$ and $i_{n-1}$ is represented by the obvious quotient map. The 
lower sequence in the above diagram is the pull-back of the upper sequence, and $i_{n-1} \circ \pi _{n-1} = \pi _n$. The transfer of a composition is the 
composition of transfers, and transfers are natural for pull-backs. Thus
$[\tau _n \circ i_{n-1}]$ is the transfer of the lower sequence composed
with the fold map
$\bigcoprod ^p _1 E \to E$. This comes out to be $p[\tau_{n-1}]$.
\qed

\vskip .2in
The Ivanov-Harer stability theorems imply that the
 map of plus constructions,
$$
B\Gamma _{g(n) -1, 2}^+ \lra  B\Gamma _{g(n)}^+
$$
is $b_n$-connected with $b_n = \lbrack \frac{ g(n) -1} {2}\rbrack$. 
Thus
$$
[ BC_{p^n} ^{(b_n)}, B\Gamma _{g(n)-1, 2}^+] \simeq [BC_{p^n}^{(b_n)},
B\Gamma _{g(n)}^+]
$$
where the superscript indicates the $b_n$-skeleton.
Let $\mu_{n,2}: BC_{p^n}^{(b_n)} \to B\Gamma ^+_{g(n)-1,2}$ correspond to
$\mu_n$ so that
$$
T\circ \omega \circ  \alpha _{n, 2} \circ \mu _{n,2}
\simeq  (1+k\psi^{-k}, (1-k) \tilde \tau_n)
\tag 5.4
$$
in $[BC_{p^n}^{(b_n)}, \Omega ^\infty _0 S^\infty (\Bbb CP ^\infty _+)]$.
Replacing $\mu_{n,2}$ by the composition
$\mu_{n+l,2} \circ i_{n+l-1} \circ \dots \circ i_{n}$ on the left side
of (5.4), by Lemma 5.3, 
the second component on the
right hand side has to be replaced by
$  p^l (1-k) \tilde \tau_n$.
This second component
is an element in $[ BC_{p^n} ^{(b_n)}, 
(\Omega_0 ^\infty S^\infty)^\wedge _p]$ which is a 
finite $p$-group. Hence, for some $l$ large enough, $p^l (1-k)  \tilde \tau_n$
is homotopic to zero. 

Consider the subset $G_n$ of $[BC_{p^n}^{(b_n)}, (B\Gamma ^+_{\infty,2})^\wedge
_p]$ of elements  that satisfy (5.4) with the second component on the
right hand side actually zero. As we have argued it is non-empty. Furthermore,
any two elements in $G_n$  differ by a map into  the fibre 
of $\omega$. The set $[ BC_{p^n} ^{(b_n)}, \Omega_0 ^\infty
S^{\infty -2}]$  is however finite and hence $G_n$ is finite. By
Tychonov's theorem, the inverse limit  of the $G_n$ is therefore non-empty. 
We pick $\tilde \mu_p \in \underset \gets \to \lim G_n$. 
Since $\colim _{n\to \infty} BC_{p ^n} ^{(b_n)} = 
BC_{p^\infty}$ has $p$-adic completion homotopy equivalent 
 to $(\Bbb C P^\infty)^\wedge
_p$, the map $\tilde \mu_p$ extends to a map
$$
\tilde  \mu_p: \Omega ^\infty S^\infty (\Bbb CP^\infty) 
\lra (B\Gamma ^+_{\infty,2} )^\wedge_p
$$
and
$$
T\circ \omega \circ \alpha \circ \mu_p \simeq (1+ k\psi ^{_k},0).
$$
This proves that the first column in the matrix described
in Theorem 2.2 is $(1+k \psi ^{-k}, 0)$.
The second column corresponds to the homotopy class of
$$
\Omega ^\infty S^\infty \overset \mu \to \lra \Bbb Z \times 
B\Gamma ^+_{\infty ,2} \overset \alpha \to \lra \Omega ^\infty 
\Bbb CP^\infty _{-1} \overset \omega \to \lra 
\Omega ^\infty S^\infty (\Bbb CP^\infty _+).
\tag 5.5
$$
An infinite loop map with source $\Omega ^\infty S^\infty $ is determined 
by its restriction to
$S^0=\{ -1, +1\} \hookrightarrow \Omega ^\infty S^\infty$. The $\mu$
above maps the base point $+1$ of $S^0$ into $(0,*)$ and the non-base point 
into $(1,*)$. The composition (5.5) maps the base point into the base point of
$\Omega ^\infty _0 S^\infty (\Bbb CP^\infty)$ and the non-base point into the
base point of $\Omega _{-2}^\infty S^\infty (\Bbb CP^\infty_+)$. This shows
 that the second column of the matrix in Theorem 2.2 is $(0,-2)$ as claimed.
\qed

\vskip .2in
\subheading{5.2. The referee's proof}
By Theorem 2.4, we may identify $\Bbb Z\times B\Gamma ^+_\infty$ with
$\Omega ^\infty \Bbb CP^\infty _{-1}$ via the map $\alpha$. Hence it suffices 
to show that the self-map 
$$
(1+k\psi ^{-k}, -2): 
\Omega ^\infty S^\infty (\Bbb CP^\infty) \times \Omega ^\infty S^\infty 
\lra
\Omega ^\infty S^\infty (\Bbb CP^\infty) \times \Omega ^\infty S^\infty 
$$ 
lifts through
$$\omega : \Omega ^\infty 
\Bbb CP^\infty _{-1} \to \Omega ^\infty S^\infty (\Bbb CP^\infty _+)
\simeq 
\Omega ^\infty S^\infty (\Bbb CP^\infty) \times \Omega ^\infty S^\infty 
$$ 
after $p$-completion. By the discussion in section 2.1,
the factor of $-2$ on the summand  $\Omega ^\infty S^\infty $ can be lifted 
by sending the non-base point of $S^0$ to the base point of the 1-component in
$\Omega ^\infty \Bbb CP^\infty _{-1}$, and extending this map to the uniquely
determined
infinite loop space map from $\Omega ^\infty S^\infty$. Hence, 
by the fibration sequence (2.3), we are left to show that after
$p$-completion the composite map
$$
\Omega ^\infty S^\infty (\Bbb CP^\infty)
\overset {1+k\psi ^k} \to \lra 
\Omega ^\infty S^\infty (\Bbb CP^\infty)
\overset \partial  \to \lra \Omega ^\infty S^{\infty -1}.
\tag 5.6
$$
is null-homotopic.

\vskip .1in
We recall  that $\partial$ is the $S^1$-transfer map, cf. [R].
Stefan Stolz gave the following, alternative description of the $S^1$-transfer
map.
Let $l: \Bbb CP ^\infty \to BU$  the map that represents the reduced
canonical line bundle, 
$\beta : BU \to \Omega U$ be the Bott map, and
$J: U \to \Omega ^\infty S^\infty$ be the $J$-homomorphism that associates
to a unitary transformation the induced map of spheres.

\vskip .2in
\proclaim{ Lemma 5.4 [St]}
The map $\partial$ is homotopic to the infinite loop space map induced by 
the composite map
$$
\Bbb CP^\infty \overset { l} \to \lra BU \overset {\beta} \to \lra 
\Omega  U \overset {\Omega J } \to \lra \Omega ^\infty S ^{\infty -1}.
$$
\endproclaim

The maps $\psi^{-k}$ give rise to  the Adams operations $\psi ^{-k}:
BU \to BU$. These operations are stable in the sense that 
they satisfy the idenity
$$
\Omega ^2 (\psi ^{-k}) \circ \beta \simeq \beta \circ (-k\psi ^{-k}).
$$
Here 
we identify $U \simeq \Omega BU$. 
Hence, we have a commutative diagram
$$
\CD
\Bbb CP^\infty @> 1 + k \psi ^{-k} >> \Omega ^\infty S^\infty (\Bbb CP^\infty)
	@> \partial >> \Omega ^\infty  S^{\infty -1}\\
@V l VV	@V lVV	@| 	\\
BU @> 1+ k \psi ^{-k} >> BU @.	\\
@V \beta VV	@V\beta VV	@|	\\
\Omega ^2 BU @> \Omega ^2 (1- \psi ^{-k})>>	\Omega ^2 BU
	@> \Omega ^2 B J >> \Omega ^\infty S^{\infty -1}.
\endCD
$$
By the afirmed Adams conjecture, 
after inverting $-k$ and hence after $p$-completion,
$BJ \circ (1 - \psi ^{-k})$ is null-homotopic.
Thus the bottom row  defines a null-homotopic map after $p$-completion which
proves (5.6).
\qed

\vskip .3in
\subheading {5.3. Proof of Theorem 1.5}.  It follows from [MS,
theorem 7.8] that for odd primes $p$ there is a splitting of
$p$-localized spaces
$$
(\Omega_0^\infty\C P^\infty_{-1})_{(p)} \simeq T_{(p)} \times BU_{(p)},
$$
where $T_{(p)}$ has torsion homotopy and homology groups.  Therefore
we have
$$
H^*_{free}(\Omega^\infty_0\C P^\infty_{-1};\Z_{(p)}) \simeq
      H^*_{free}(BU;\Z_{(p)}) = H^*(BU;\Z_{(p)}).
$$

\vskip .1in
For $p=2$ we use the fact
that $H_*(\Omega^\infty\C P^\infty_{-1};\F_2)$ is a
sub Hopf algebra of $H_*(\Omega ^\infty S ^\infty
\C P^\infty_+;\F_2)$ by Theorem 1.3 of [G1].
We need explicit additive generators in low dimensions.  Let $a_i \in
H_{2i}(\Omega ^\infty S ^\infty
_1 \C P^\infty_+;\F_2)$ be the non-zero element in the image of
$\C P^\infty \to \Omega ^\infty S^\infty_1\C P^\infty_+$ into the 1-component of
$\Omega ^\infty S^\infty \C P^\infty_+$ and let
$$
Q^s: H_*(\Omega ^\infty S^\infty \C P^\infty_+; \F_2) \to H_{*+s}(
\Omega ^\infty S^\infty \C P^\infty_+; \F_2)
$$
be the $s$-th  homology operation for infinite loop spaces.  
The following elements form an
$\F_2$-basis for $H_*(\Omega^\infty_0\C P^\infty_{-1};\F_2)$,
$$
\aligned
x_2 &= (Q^1 a_0)^2a_0^{-4},\\
x_3 &= (Q^3a_0)a_0^{-2} + (Q^2Q^1a_0)a_0^{-4} +
	  (Q^2a_0)(Q^1a_0)a_0^{-4} + (Q^1a_0)^3a_0^{-6},\\
x_4 &= (a_1)^2 a_0^{-2}, \quad y_4 = (Q^2a_0)^2a_0^{-4}, \quad x_2^2,
\endaligned
$$
in positive degrees less than 5.
The rational homology of $\Omega^\infty_0\C P^\infty_{-1}$ and $BU$
agree, so
$$
H^{free}_2(\Omega^\infty_0\C P^\infty_{-1}) = \Z, \, 
H^{free}_3(\Omega^\infty_0\C P^\infty_{-1}) = 0, \,
H^{free}_4(\Omega^\infty_0\C P^\infty_{-1}) = \Z\oplus\Z.
$$
We now employ the Bockstein spectral sequence.  For an infinite loop
space $X$, it is a singly graded spectral sequence of Hopf algebras
with
$$
E^1_* = H_*(X; \F_2), \quad E^\infty_* = H_*^{free}(X) \otimes \F_2.
$$
The differential $d^1$ is the standard Bockstein [Br].  The elements
listed above for $X = \Omega^\infty_0\C P^\infty_{-1}$ all survive to
$E^2_*$.  The general formula for $d^2$ applied to an even square [M,
Proposition 1.5] gives
$$
d^2 y_4 = x_3 + (Q^3a_0)a_0^{-2} + (Q^1 a_0)^3a_0^{-6} = x_3
$$
in $E^2_*(\Omega ^\infty S^\infty
_0\C P^\infty_+)$, and $d^2 x_4 = d^2 x_2^2 = 0$.  On the
other hand, the $E^\infty_*$-term is $\F_2$, 0, and $\F_2 \oplus \F_2$
in degrees 2, 3 and 4, respectively.  It follows that
$$
E^\infty_2(\Omega^\infty_0\C P^\infty_{-1}) = \F_2\langle
x_2\rangle,\quad
E^\infty_4(\Omega^\infty_0\C P^\infty_{-1}) = \F_2\langle x_4\rangle
\oplus \F_2\langle x_2^2\rangle.
$$
Both $x_2^2$ and $x_4$ are primitive elements (since that is the case
in $\Omega ^\infty S^\infty
_0\C P^\infty_+$), so the reduced diagonal from $E^\infty_4$ to
$E_2^\infty \otimes E_2^\infty$ is zero.  Dually the squaring map is
zero on $H^2_{free}(\Omega^\infty_0\C P^\infty_{-1}) \otimes \F_2$, so
$H^*_{free}(\Omega^\infty_0\C P^\infty_{-1})$ is not polynomial.\qed

\vskip .2in
\subheading{Remark}. The non-zero differential $d^2: E^2_4 \to 
E^2_3$
proves that the two-torsion in $H_3(B\Gamma_\infty)$ is $\Bbb Z/4\Bbb Z$.

\vskip .4in
\subheading {References}

\vskip .2in
\noindent
[A1] J.F. Adams,
Lectures on generalized cohomology, category theory, homology
theory and their applications III, LNM 99, Springer (1969).

\vskip .1in
\noindent
[A2] J.F. Adams, {\it On the groups $J(X)$ II}, Topology {\bf 3} (1965),
137--171.

\noindent
[Ak] T. Akita, {\it Nilpotency and triviality of mod $p$ Morita-Mumford
classes of mapping class groups of surfaces}, Nagoya Math. J. {\bf 165 }
(2002), 1--22.

\vskip .1in
\noindent
[BS] J.C. Becker, R.E. Schultz, {\it Axioms for transfers and traces}, Math. Z. {\bf 227} (1988), 583-605.

\vskip .1in
\noindent
[B] J. Berrick, An approach to algebraic $K$-theory, Research Notes in 
Mathematics, 56, Pitman (Advanced Publishing Program), Boston, 1982.

\vskip .1in
\noindent
[BK] A.K. Bousfield, D.M. Kan, Homotopy Limits, Completions and Localizations,
Springer Verlag LNM {\bf 304}, 1972.

\vskip .1in
\noindent
[Br] W. Browder, {\it Torsion in $H$-spaces}, Ann. of Math. {\bf 74} (1961),
24-51.

\vskip .1in
\noindent
[C] M. Crabb, {\it $\Bbb Z /2$-Homotopy Theory}, LMS Lecture Note Series 
{\bf 44}, CUP, 1980.

\vskip .1in
\noindent
[EE] C.J.  Earle, J. Eells, {\it A fibre bundle description of
Teichm\"uller theory}, J. Diff. Geom. {\bf 3} (1969), 19--43.

\vskip .1in
\noindent
[ES] C.J. Earle, A. Schatz, {\it Teichm\"uller theory for surfaces
with boundary}, J. Diff. Geom. {\bf 4} (1970), 169--185.

\vskip .1in
\noindent
[F] O. Forster: Lectures on Riemann surfaces, Graduate Texts in
Mathematics, 81. Springer-Verlag, New York-Berlin, 1991.

\vskip .1in
\noindent
[G1] S. Galatius, {\it  Mod $p$ homology of the stable mapping class group},
Topology {\bf 43} (2004), 1105-1132.

\vskip .1in
\noindent
[G2] S. Galatius, {\it Secondary characteristic classes of surface bundles},
preprint, math.AT/0402226.

\vskip .1in
\noindent
[H1] J.L. Harer, {\it The second homology group of the mapping class group
of an orientable surface }, Invent. Math. {\bf 72} (1983), 221--239.

\vskip .1in
\noindent
[H2] J.L. Harer, {\it Stability of the homology of the mapping class groups of
orientable surfaces}, Annals Math. {\bf 121} (1985), 215--249.

\vskip .1in
\noindent
[I] N.V. Ivanov, {\it Stabilization of the homology of Teichm\"uller modular
groups},
Original: Algebra i Analiz {\bf1} (1989), 110--126; Translated:
Lenin\-grad Math. J. {\bf 1} (1990), 675--691.

\vskip .1in
\noindent
[LMS] G. Lewis, J.P. May, M. Steinberger, {\it Equivariant stable 
homotopy theory}, Springer LNM {\bf 1213}, 1986.

\vskip .1in
\noindent
[M] I. Madsen, {\it  Higher Torsion in $SG$ and $BSG$}, Math. Z.
{\bf 143} (1975), 55-80.

\vskip .1in
\noindent
[MS] I. Madsen, C. Schlichtkrull, {\it  The circle transfer and $K$-theory},
AMS Contemporary Math. {\bf 258} (2000), 307--328.

\vskip .1in
\noindent
[MT] I. Madsen, U. Tillmann, {\it The stable mapping class group
and $Q(\Bbb CP^\infty)$}, Invent. Math. {\bf 145} (2001), 509--544.

\vskip .1in
\noindent
[MSt] J. Milnor, J. Stasheff, Characteristic Classes,  Study 76,
Princeton University 
Press 1974.

\vskip .1in
\noindent
[MW] I. Madsen, M. Weiss, {\it The stable moduli space of Riemann surfaces:
Mumford's conjecture}, math.AT/0212321.

\vskip .1in
\noindent
[Mi] E.Y. Miller, {\it The homology of the mapping class group}, J. Diff. Geom.
{\bf 24} (1986), 1--14.

\vskip .1in
\noindent
[Mo] S. Morita, {\it Characteristic classes of surface bundles}, Invent. Math.
{\bf 90} (1987), 551-577.

\vskip .1in
\noindent
[Mu] D. Mumford, {\it Towards an enumerative geometry of the moduli space of
curves}, in Arithmetic  and Geometry, M. Artin and J. Tate, editors, Progr.
Math. {\bf 36} Birkhauser (1983), 271--328.

\vskip .1in
\noindent
[R] D. Ravenel, {\it The Segal conjecture for cyclic groups and its
consequences}, Amer. J. Math. {\bf 106} (1984), 415--446.

\vskip .1in
\noindent
[S] G. Segal, {\it The stable homotopy of complex projective space},
Quart. J. Math. Oxford {\bf 24} (1973), 1--5.

\vskip .1in
\noindent
[St] S. Stolz, {\it Relationships between transfer and $J$-homomorphisms},
Bonner Mathematische Zeitschriften {\bf 125}, Universit\"at Bonn, 1980.

\vskip .1in
\noindent
[T] U. Tillmann, {\it On the homotopy of the stable mapping class group},
Invent. Math. {\bf 130} (1997), 257--275.

\parindent =0in
\vskip .3in

Soren Galatius and Ib Madsen

Matematisk Institut

Aarhus Universitet

8000 Aarhus C

Denmark


\vskip .2in
Ulrike Tillmann

Mathematical Institute

24-29 St. Giles Street

Oxford OX1 3LB

UK


\enddocument